\documentclass[12pt]{amsart}
\usepackage{amscd, amssymb}

\newtheorem{pr}{Proposition}
\newtheorem{lm}{Lemma}

\newcommand{\proj}{\mathbf P}

\newcommand{\barr}{\overline}
\newcommand{\rarr}{\rightarrow}
\newcommand{\oh}{{\mathcal{O}}}
\newcommand{\com}{\mathbb{C}}
\newcommand{\G}{\mathbf{G}}
\newcommand{\T}{\mathbf{T}}
\newcommand{\Q}{\mathbb{Q}}

\newcommand{\Z}{\mathbb{Z}}

\newcommand{\Xvir}{[X]^{\it{vir}}}
\newcommand{\Xivir}{[X_i]^{\it{vir}}}
\newcommand{\Nvir}{N^{\it{vir}}}
\newcommand{\Dvir}{D^{\it{vir}}}
\newcommand{\Divir}{D^{\it{vir}}_i}
\newcommand{\A}{\mathbf{A}}

\newcommand{\eqq}{\stackrel{\sim}{=}}

\newcommand{\sumo}{\oplus}

\newcommand{\bpf}{\noindent {\em Proof.} }
\newcommand{\epf}{\qed \vspace{+10pt}}

\newcommand{\eref}{e_{\rm{ref}}}

\newcommand{\Ext}{{\rm Ext}}
\newcommand{\Def}{{\rm Def}}
\newcommand{\val}{{\rm val}}
\newcommand{\Aut}{{\rm Aut}}
\newcommand{\C}{\mathbb{C}}
\newcommand{\rk}{\mbox{rk}}

\newcommand{\M}{{\barr M}}
\newcommand{\Mgn}{\M_{g,n}}
\newcommand{\Mgnp}{\Mgn (\proj^r,d)}

\begin{document}
\title{Localization of Virtual Classes}
\author{T. Graber and R. Pandharipande}
\date{28 August 1997}
\maketitle

%\begin{center}
%{\it Preliminary version}
%\end{center}

\pagestyle{plain}
\setcounter{section}{-1}

\section{\bf{Introduction}}
We prove  a localization formula for the
virtual fundamental class in the general context of
$\com^*$-equivariant perfect obstruction theories.
Let $X$ be an algebraic scheme with a $\com^*$-action
and a $\com^*$-equivariant perfect obstruction theory.
The virtual fundamental class $\Xvir$ in the expected
equivariant Chow group $A_*^{\com^*}(X)$ may be constructed
by the methods of Li-Tian [LT] and Behrend-Fantechi [B], [BF].
Each  connected component $X_i$ 
of the fixed point scheme carries an associated 
$\com^*$-fixed perfect obstruction theory. A 
virtual fundamental class in $A_*(X_i)$ is thus determined.
The virtual normal bundle to $X_i$ is
obtained from the moving part of the virtual tangent space
determined by the obstruction theory. 
The localization formula is then:
\begin{equation}
\label{exloc} \Xvir =
\iota_* \sum  \frac{\Xivir}{e(N^{\it{vir}}_i)}
\end{equation}
in $A_*^{\com^*}(X) \otimes \Q[t,\frac{1}{t}]$ where
$t$ is the generator of the equivariant ring of $\com^*$.
This localization formula is the main result of the paper.
The proof 
requires an additional hypothesis on $X$:
the existence of a $\com^*$-equivariant embedding in a
nonsingular variety $Y$.

In case $X$ is nonsingular with the trivial perfect
obstruction theory, equation (\ref{exloc}) 
reduces immediately to the standard localization formula
[Bo], [AB]. Originally, this localization was
proven in equivariant cohomology. Algebraic localization in
equivariant Chow theory has recently been established in
[EG2]. The point of view of our paper is entirely algebraic.

The definitions and constructions related to the virtual
localization formula (\ref{exloc}) are
discussed in Section \ref{locfor}. The simplest example
of a $\com^*$-equivariant perfect obstruction theory
is given by the following situation.
Let $Y$ be a nonsingular algebraic variety with a
$\com^*$-action. Let $V$ be a $\com^*$-equivariant 
bundle on $Y$. Let $v\in H^0(Y, V)^{\com^*}$ be a $\com^*$-fixed
section. Let $X$ be the scheme-theoretic zero locus of $v$.
$X$ is naturally endowed with an equivariant 
perfect obstruction theory which yields the refined
Euler class (top Chern class)  of $V$ as the virtual fundamental class.
The localization formula in this basic setting
is proven in Section \ref{localcase}. The method is
to deduce (\ref{exloc}) for $X$ from the
known localization formula for the nonsingular
variety $Y$.

The proof of (\ref{exloc}) for general
$\com^*$-equivariant perfect obstruction theories on an
algebraic scheme $X$
proceeds in a similar manner. Again,
formula (\ref{exloc}) is deduced from the
ambient localization formula for $Y$.
The argument here is more
subtle: explicit 
manipulation of cones and a rational equivalence
due to Vistoli [V] are necessary.
This proof is given in Section \ref{gencase}.

There are
two immediate applications of the virtual 
localization formula. First, a local complete
intersection scheme is endowed with a
canonical perfect obstruction theory obtained
from the cotangent complex. A localization
formula is thus obtained for these singular
schemes (at least when equivariant
embeddings in nonsingular varieties exist).
Second, the proper Deligne-Mumford moduli stack
$\overline{M}_{g,n}(V, \beta)$ of stable maps to
a nonsingular projective variety $V$ is equipped
with a canonical perfect obstruction theory. 
If $V$ has a $\com^*$-action, then a natural
$\com^*$-action on $\overline{M}_{g,n}(V, \beta)$
is defined  by translation of the map.
An equivariant perfect obstruction theory 
on $\overline{M}_{g,n}(V, \beta)$ can be obtained.
Moreover,  
$\overline{M}_{g,n}(V, \beta)$ admits an equivariant embedding
in a nonsingular Deligne-Mumford
stack.
As a result, the virtual localization formula
holds for $\overline{M}_{g,n}(V, \beta)$.

In the last two sections of the paper,  consequences of the
localization formula in Gromov-Witten theory are explored.
In Section \ref{projj}, an explicit graph summation formula
for the Gromov-Witten invariants (for all genera) of
$\proj^r$ is presented via localization on the
moduli space of maps $\overline{M}_{g,n}(\proj^r,d)$. 
The invariants are expressed as a sum over graphs 
corresponding to the fixed point loci. For each graph,
the summand is a product over vertex terms.
The vertex terms  
are integrals over associated spaces 
$\overline{M}_{g',n'}$ of the Chern classes of the
cotangent line bundles and the Hodge bundle.
All these integrals may be calculated from Witten's conjectures
(Kontsevich's theorem) by a method due to Faber [Fa].
Similar graph sum formulas exist for Gromov-Witten
invariants (and their descendents) of all the
compact algebraic homogeneous spaces $\G/\mathbf{P}$ via the
action of the maximal torus $\T\subset \G$.

The positive
degree
Gromov-Witten invariants of $\proj^2$ coincide with
enumerative geometry: they count the number $N^{g}_d$
of genus $g$, degree $d$, nodal plane curves passing 
through $3d+g-1$ points in $\proj^2$.
Localization presents a solution of this enumerative
geometry problem via integrals of tautological
classes over the  moduli space of pointed curves.
The numbers $N^{g}_d$ have been computed via 
more classical degeneration methods in [R1], [CH].
The character of the solutions in [R1], [CH] is markedly different:
it is 
by recursion over wider classes of enumerative
questions.

In Section \ref{ellip}, localization is applied to a question
suggested to us by S. Katz:
the calculation of excess integrals on the moduli spaces
$\overline{M}_{g,0}(\proj^1,d)$ that arise in the study of
Calabi-Yau 3-folds. Under suitable conditions,
the integral
\begin{equation}
\label{exxx}
\int _{[ \barr{M}_{g,0}(\proj^1,d)]^{\it{vir}}} 
c_{\rm{top}} (R^1 \pi_* \mu^* N)
\end{equation}
is the contribution to the
genus $g$ Gromov-Witten invariant of a Calabi-Yau 3-fold 
of multiple
covers of a fixed rational curve (with normal bundle
$N=\oh(-1)\oplus \oh(-1)$). In [M], the integral (\ref{exxx})
is explicitly evaluated to be $1/d^3$ 
in the genus $g=0$ case via localization
on the nonsingular stack $\barr{M}_{0,0}(\proj^1,d)$.
A trick of setting one of the $\com^*$-weights on $\proj^1$
to be 0 is used. We evaluate the excess integral
in the genus $g=1$ case
in Section \ref{ellip} via virtual localization on
$\barr{M}_{1,0}(\proj^1,d)$.
Manin's trick [M] 
and formulas for cotangent line integrals on $\overline{M}_{1,n}$
are used to handle the graph sum.
The answer obtained, $1/12d$, agrees with the physics result of
[BCOV]. 
The higher genus  integrals may be explicitly
evaluated in any given case by virtual localization and the
algorithm implemented by Faber [Fa]
to calculate the vertex integrals. The conjecture obtained
from these calculations is: for $g\geq 2$,
\begin{equation}
\label{conjj}
\int _{[ \barr{M}_{g,0}(\proj^1,d)]^{\it{vir}}} c_{\rm{top}} 
(R^1 \pi_* \mu^* N)
= \frac{|B_{2g}| \cdot d^{2g-3}}{2g\cdot (2g-2)!} =
|\chi(M_g)| \frac{d^{2g-3}}{(2g-3)!}
\end{equation}
where $B_{2g}$ is the $2g^{th}$ Bernoulli number and
$\chi(M_g)= B_{2g}/2g(2g-2)$ is the orbifold
Euler characteristic of $M_g$.
This conjecture was made jointly with C. Faber.
We have not yet been able to evaluate the
graph sums uniformly to establish (\ref{conjj}).

The localization formula and graph sum formulas were 
first introduced in the context of stable maps  by Kontsevich
in [K] following related work of Ellingsrud and Str\o mme [ES].
Kontsevich studied the convex genus 0 case
 where the moduli spaces are nonsingular Deligne-Mumford stacks.
Many ideas about the virtual fundamental class
and localization described here are implicit in [K].
In particular, the  higher genus formulas of Section
\ref{projj} are identical to the genus 0 formulas of [K]
except for the new Hodge bundle terms. However, the higher
genus map spaces are in general nonreduced, reducible, and
singular, so the virtual localization formula (\ref{exloc})
is essential.  
 Givental has stated
a localization axiom for genus 0 Gromov-Witten
invariants of toric varieties in [G] which follows from
(\ref{exloc}). Localization formulas are used in [G]
to prove predictions of mirror symmetry in the
case of Calabi-Yau complete intersections in
toric varieties.

The authors thank P. Aluffi,
K. Behrend, D. Edidin,
W. Fulton, E. Getzler, L. G\"ottsche, S. Katz, A. Kresch, J. Li,
B. Seibert, M. Thaddeus, and A. Vistoli for many valuable conversations. 
A special thanks is due to C. Faber for his
computations of the vertex integrals in (\ref{conjj}) and to
B. Fantechi for her
tireless explanations of obstruction theories and
virtual classes.
The first author was supported by an NSF graduate fellowship.
The second author was partially supported by an NSF
post-doctoral fellowship.
The authors also thank the
Mittag-Leffler Institute for support.

\section{\bf{The virtual localization formula}}
\label{locfor}
Let $X$ be an algebraic scheme over $\com$.
A perfect obstruction theory consists of the following data:
\begin{enumerate}
\item[(i)] A two term complex of vector bundles 
$E^\bullet = [E^{-1} \rarr E^0]$ on X.
\item[(ii)] A morphism $\phi$ in the derived category (of
quasi-coherent sheaf complexes bounded from above)
from $E^\bullet$
to the cotangent complex $L^\bullet X$ of $X$ satisfying two properties.
\begin{enumerate}
\item[(a)] $\phi$ induces an isomorphism in cohomology in degree 0.
\item[(b)] $\phi$ induces a surjection in cohomology in degree -1.
\end{enumerate}
\end{enumerate}
The constructions
of [LT],  [BF]  give rise to a virtual fundamental class, $\Xvir$
in $A_d(X)$ where $d= \rk(E^0)-\rk(E^{-1})$
is the expected dimension. Let $E_\bullet=[ E_0 \rarr E_1]$
denote the dual complex of $E^\bullet$.

If $X$ admits a global closed embedding in a nonsingular
scheme  (or Deligne-Mumford stack) $Y$, one can give a relatively
straightforward construction of the 
virtual class as follows.
In this situation, the two term cut-off of the
cotangent complex can be taken to be:
$$L^\bullet X = [I/I^2 \rarr \Omega_Y]$$
where $I$ is the ideal sheaf of $X$  in $Y$.
Since only this cut-off will be used, the cotangent
complex will be identified with its cut-off throughout this
section.
We assume for simplicity that 
\begin{equation}
\label{assum}
\phi: E^\bullet \rarr [I/I^2 \rarr \Omega_Y]
\end{equation}
is an actual map of complexes.
This hypothesis is not required for the constructions
of [LT], [BF].
However, if $X$ has enough locally frees,
such  a representative $(E^\bullet, \phi)$
may always be chosen in the derived category.

The mapping cone associated to the morphism $\phi$ of complexes
 yields an exact sequence
of sheaves:
\begin{equation}
\label{mppp}
E^{-1} \rarr E^0 \sumo I/I^2 \stackrel{\gamma}{\rarr} \Omega_Y \rarr 0.
\end{equation}
In fact, $\phi$
satisfies (a) and (b) if and only if 
 (\ref{mppp}) is exact.
We consider the associated exact sequence of abelian cones
\begin{equation}
\label{fff}
0 \rarr TY \rarr C(I/I^2) \times_X E_0 \rarr C(Q)\rarr 0
\end{equation}
where $C(Q)$ is the cone associated to the kernel $Q$ of $\gamma$. 
As $Q$ is a quotient of $E^{-1}$, 
$C(Q)$ embeds in $E_1$.  The normal cone of $X$ in $Y$, $C_{X/Y}$, 
is naturally a closed subscheme of  $C(I/I^2)$.
If we
define $D=C_{X/Y} \times_X E_0$, then $D$ is a $TY$-cone (see [BF]), and the
quotient of $D$ by $TY$ is a subcone of $C(Q)$ which we will denote
by $D^{\it{vir}}$.  The virtual fundamental class of $X$ associated
to this obstruction theory is then the refined
intersection of $D^{\it{vir}}$ with
the zero section of the vector bundle $E_1$.

Suppose $X\subset Y$ is equipped with an equivariant  $\G$-action together
with a lifting to the complex $E^\bullet$ such that
$\phi$ is a morphism in the derived category of $\G$-equivariant
sheaves (with respect to the natural $\G$-action on
$L^\bullet X$). 
The above construction then
yields an equivariant virtual fundamental class in the 
equivariant Chow group $A^\G_d(X)$ since the cones used are invariant.
In fact, to define the equivariant virtual class,
global equivariant embeddings
are not necessary.

We now assume the group $\G$ is
the torus
$\com^*$.  We expect to be able to reduce integrals over $\Xvir$
to integrals over the fixed point set. 
Let $X^f$ be the maximal $\com^*$-fixed closed subscheme
of $X$. $X^f$ is the natural scheme theoretic fixed point
locus. 
If $X={\rm Spec}(A)$, then the ideal of $X^f$ is generated by
the $\com^*$-eigenfunctions
with nontrivial characters.
For nonsingular $Y$, $Y^f$ is the nonsingular
set theoretic fixed point locus [I].
For $X\subset Y$, the relation
$X^f= X\cap Y^f$ holds. We let $Y_f= \bigcup Y_i$
be the decomposition into irreducible components.
Let $X_i= X\cap Y_i$. $X_i$ is possibly reducible.

Let $S$ be a coherent sheaf on a fixed component $X_i$ with a
$\com^*$-action. $S$ decomposes as direct sum,
$$S = \bigoplus_{k\in \mathbb{Z}} S^{k},$$
of $\com^*$-eigensheaves of $\oh_{X_i}$-modules.
If $S$ is locally free, each summand is also locally
free. We denote the fixed subsheaf $S^{0}$ by $S^f$ and the
moving subsheaf
$\oplus_{k\neq 0} S^k$ by $S^{m}$. 

There
is a natural isomorphism $\Omega_Y |_{Y_i}^f= 
\Omega_{Y_i}$ [I]. It is easy to then deduce:
$$\Omega_X |_{X_i}^f = \Omega_{X_i}$$
from the equality $X_i = X \cap Y_i$ or the universal
property of K\"ahler differentials.

Let $E^{\bullet}_i$ denote the restriction of $E^\bullet$
to $X_i$. Let $E^{\bullet,f}_i$ denote the
fixed part of the complex $E^{\bullet}_i$. $E^{\bullet,f}_i$
is a two term complex of bundles.
There exists a canonical map,
$$\psi_i: E^{\bullet,f}_i  \rarr L^\bullet X_i,$$
determined by the following
construction.
Let $\phi_i: E^\bullet_i \rarr L^\bullet X |_{X_i}$ be
the pull-back of $\phi$, and let 
$\phi_i^f: E_i^{\bullet,f} \rarr L^\bullet X|_{X_i}^f$ be
the associated fixed map.
Similarly let $\delta_i: L^\bullet X|_{X_i} \rarr
L^\bullet X_i$
be the canonical morphism, and let $\delta_i^f$ be
the associated fixed map. Then,
$\psi_i= \delta_i^f \circ \phi_i^f$.

\begin{pr}
The map 
$\psi_i: E^{\bullet,f}_i  \rarr L^\bullet X_i$ is a 
canonical perfect obstruction theory on $X_i$.
\end{pr}

\bpf
To show $\psi_i$ satisfies properties (a) and
(b), it suffices to show both maps $\phi_i^f$ and
$\delta_i^f$ satisfy these properties. 

A map of complexes $\nu: A^\bullet \rarr B^\bullet$
satisfies (a) and (b) if and only if the sequence
$$ A^{-1} \oplus B^{-2} \rarr A^0 \oplus B^{-1} \rarr B^0 \rarr 0$$
 is exact. Since tensor product is
right exact, the joint validity of (a) and (b) is  preserved
under pull-back. As $\phi$ is a perfect
obstruction theory, $\phi_i$ satisfies (a) and (b).
The fixed map $\phi_i^f$
also satisfies properties (a) and (b) since taking invariants
is exact.

The cotangent
complex of $X_i$ can be represented by the embedding $X_i \subset Y_i$:
$$L^\bullet X_i =[ I_{X_i/Y_i}/ I^2_{X_i/Y_i} \rarr \Omega_{Y_i} 
|_{X_i}].$$
The zeroth cohomology of $L^\bullet X|_{X_i}^f$ is
$\Omega_{X}|_{X_i}^f = \Omega_{X_i}.$
Thus, $\delta_i^f$ satisfies property (a).

Property (b) for $\delta^f_i$ will now be established.
The map $\delta^f_i$ is represented by the natural diagram:
\begin{equation*}
\begin{CD}
 I_{X/Y}/I^2_{X/Y} | _{X_i}^f 
@>>>  \Omega_{Y} |_{X_i}^f \\
@V{d^{-1}}VV   @V{d^0}VV \\
 I_{X_i/Y_i}/I^2_{X_i/Y_i} @>>> \Omega_{Y_i}|_{X_i}
\end{CD}
\end{equation*}
Since $X_i= X\cap Y_i$, the
map $$I_{X/Y}/I^2_{X/Y}|_{X_i} \rarr I_{X_i/Y_i}/I^2_{X_i/Y_i}$$
is surjective. 
Hence, $d^{-1}$ is surjective. As
$d^0$ is an isomorphism,
$\delta_i^f$ 
is surjective on cohomology in degree $-1$.
\epf

\noindent The virtual structure on $X_i$ is defined to be the
one induced by the perfect obstruction theory
$\psi_i: E_i^{\bullet,f} \rarr L^\bullet X_i$.

We define the virtual normal bundle, $\Nvir_i$ to $X_i$ to be 
the moving part of 
$E_{\bullet,i}$.  
Note that $E_{\bullet,i}$ is a complex, and not a single bundle.  
Also note that in the non-virtual case, when the complex has just
one term, this coincides with the usual normal bundle.
Define the Euler class (top
Chern class) of a two term complex $[B_0 \rarr B_1]$ to be
the ratio of the Euler classes of the two bundles:
$e(B_0)/e(B_1)$. We arrive at
the following natural formulation of the virtual Bott residue formula
for the Euler class of a bundle $A$ of rank equal to the virtual dimension
of $X$:
\begin{equation}
\label{boott}
\int_{\Xvir} e(A) = \sum \int_{\Xivir} \frac{e(A_i)}{e(\Nvir_i)}
\end{equation}
where the Euler classes on the right hand side are
equivariant classes. Since  $\Nvir_i$ is
a complex of bundles with nonzero $\com^*$-weights, the
Euler class $e(\Nvir_i)$ is invertible in the localized
ring 
$$A^{\com^*}_*(X)_{{t}}= A^{\com^*}_*(X) \otimes_{\Q[t]}\Q[t,
\frac{1}{t}].$$
Chow groups will always be taken with $\Q$-coefficients.
As in the case of the standard Bott residue formula, equation (\ref{boott})
should
be a consequence of a localization formula in equivariant Chow groups.
On a nonsingular variety $Y$, the fundamental
result which immediately implies the residue formula is: 
$$[Y]= \iota_* \sum \frac{[Y_i]} {e(N_i)}$$
in $A^{\com^*}_* (Y) _t.$
The obvious generalization to the virtual setting which would just
as readily imply the virtual residue formula is:
\begin{equation}
\label{vbott}
\Xvir = \iota_* \sum \frac{\Xivir} {e(\Nvir_i)}.
\end{equation}

It is worth remarking that in the case of most interest to us,
the moduli space of maps to projective space, the right side
of (\ref{vbott}) is directly accessible.  In this case,
three special properties hold. First, the fixed loci
$X_i$ for a general $\com^*$-action
have been identified by Kontsevich in [K]: they
are indexed by graphs and are essentially products of
Deligne-Mumford moduli spaces of pointed curves.
Second,
$\Xivir=[X_i]$. Finally,  $e(\Nvir_i)$ is
expressible in terms of tautological classes on $X_i$
via the deformation theory of curves and maps.
Thus (\ref{vbott})
provides a concrete way to calculate virtual integrals on moduli
spaces where it seems quite difficult to directly evaluate the
virtual fundamental class.

\section{\bf{Proof in the basic case}}
\label{localcase}
As a first motivational step, we prove the virtual
localization formula (\ref{exloc}) in the 
following situation.
Let $Y$ be 
a nonsingular  variety equipped with a $\com^*$-action, 
a $\com^*$-equivariant bundle $V$, and an invariant 
section $v$ of $V$. The zero scheme $X$ of $v$ carries
a natural equivariant perfect obstruction theory.
The two term complex of bundles  on $X$,
$$E^\bullet =[V^\vee \rarr \Omega_Y],$$
is obtained from the the section $v$.
The required morphism to the cotangent complex
$L^\bullet X =[ I/I^2 \rarr \Omega_Y]$ is obtained from the
natural map
$V^\vee \rarr I/I^2$ on $X$.
The definitions show 
the virtual fundamental class in this case is just the refined Euler
class of $V$.  That is, if we consider the graph of the section,
and take its refined intersection product with the zero section, we 
get a Chow homology class supported on the zero locus $X$.  
The definitions of the virtual fundamental class for
general spaces are specifically designed to recover this refined
Euler class from the local data of the two term complex, and
to generalize this class in cases where such a geometric realization of 
the deformation complex does not necessarily exist.

In this basic situation, we can express all of the virtual objects in
the localization formula in terms of familiar data on $Y$.  
As in Section \ref{locfor},
we denote the components of the 
fixed locus of $Y$ by $Y_i$. 
$V$ splits into eigenbundles on $Y_i$.
Since $v$ is a $\com^*$-invariant section, it is necessarily a section of
the weight 0 bundle $V_i^f$.  
$Y_i$ is a smooth manifold with a vector bundle and a section
which vanishes exactly on $X_i=X\cap Y_i$. 
The associated $\com^*$-fixed
obstruction theory defined in Section \ref{locfor},
$$[(V_i^{f}) ^\vee \rarr \Omega_{Y_i}],$$ is exactly 
the perfect obstruction theory obtained from the pair $V_i^f$ and
$v\in H^0(Y_i, V_i^f)$. Note the maps to the cotangent
complex must be checked to agree.
It follows that
the virtual fundamental class of $X_i$ is the same as the
refined Euler class of $V_i^f$ on $Y_i$.

The virtual normal bundle is
by definition the moving part of the complex $[TY \rarr V]$. 
The moving part of $TY$ is just the normal bundle to $Y_i$. Hence
$\Nvir_i$ is  
the  complex 
$[N_{Y_i/Y} \rarr V_i^m]$.  By the definition of Euler class of a 
complex, we obtain: $$e(\Nvir_i)=\frac{e(N_{Y_i/Y})} { e(V_i^m)}.$$  
After substituting this expression for $e(\Nvir_i)$
into the virtual localization formula (\ref{exloc}), 
we see the equality we want to prove in $A_*^{\com^*}(X)_t$ is:
\begin{equation}
\label{goall}
\eref(V) = \iota_* \sum \frac{\eref(V_i^f) \cap e(V_i^m)} { e(N_{Y_i/Y})}
\end{equation}
where $\eref(V)$ is the refined Euler class of $V$ as
a Chow homology class on $X$.
We know by the
localization formula on $Y$:
$$[Y] = \iota_* \sum \frac{[Y_i]}{e(N_{Y_i/Y})}. $$
Intersecting both sides with $\eref(V)$ yields:
$$\eref(V) = \iota_* \sum \frac{\eref(V)\cap [Y_i]}{ e(N_{Y_i/Y})}.$$
Since taking refined Euler class commutes with pullback, the numerators
on the right hand side are just the refined Euler classes of $V_i$.  
On each component,
we have the splitting
 $V_i=V_i^f \oplus V_i^m$. Since the section lives
entirely in $V_i^f$, it follows that $\eref(V_i)=
\eref(V_i^f) \cap e(V_i^m)$.
Formula (\ref{goall}) is thus obtained. The proof of (\ref{exloc}) in the
basic case is complete.

\section{\bf{Proof in the general case}}
\label{gencase}
In this section, we prove the virtual localization formula for an
arbitrary scheme $X$ which admits an equivariant embedding
in a nonsingular scheme $Y$.
Recall from the construction of the virtual class in Section \ref{locfor}, 
the two cones
$D$ and $\Dvir$ satisfy:
\begin{equation}
\label{seq1}
0 \rarr TY \rarr D \rarr \Dvir \rarr 0
\end{equation}
\begin{equation}
\label{seq2}
D=  C_{X/Y} \times  E_0.
\end{equation}
$\Dvir$ is a embedded as a closed 
subcone of $E_1$. The 
virtual class is defined by  $\Xvir = s_{E_1}^* [\Dvir]$.
Alternatively, there is a fiber square:
\begin{equation}
\label{sqq}
\begin{CD}
TY @>>> D \\
@VVV   @VVV \\
X @>{0_{E_1}}>> E_1
\end{CD}
\end{equation}
where the bottom map is the zero section.
Then, $\Xvir= s_{TY}^* 0^{!}_{E_1} [D].$

Let $X_i= X\cap Y_i$ be defined as in Section \ref{localcase}.  
$X_i$ is a union of connected components. 
$\com^*$-fixed analogues of (\ref{seq1}) and (\ref{seq2})
hold for the embeddings
$X_i \subset Y_i$:
$$0 \rarr TY_i \rarr D_i \rarr \Divir \rarr 0,$$
$$D_i= C_{X_i/Y_i} \times E_0^f.$$
$\Divir$ is a embedded as a closed 
subcone of $E^f_1$, and $\Xivir= s_{E_1^f}^*(\Divir)$.
Since $X_i$ is possibly disconnected,
it should be noted that the ranks of the bundles $E^f_{0,i}$ and
$E^f_{1,i}$ may vary on the connected components.
The Euler classes of these bundles on $X_i$ are
taken with respect to their ranks on each component.
For notational convenience, the restriction subscript $i$
will be dropped. Similarly, the pull-backs of
$TY$ and $TY_i$ to $X_i$ will be denoted by $TY$ and $TY_i$.

The virtual localization formula for $X$ will be deduced from
localization for $Y$.  We start with the equality
$$[Y] = \iota_* \sum \frac{[Y_i]}{e(TY^m)}$$
in $A^{\com^*}_*(Y)_t$.
The refined intersection product with $\Xvir$ yields:
$$\Xvir = \iota_* \sum \frac{\Xvir \cdot [Y_i]}{e(TY^m)}$$
in $A^{\com^*}_*(X)_t$.
Comparing this equation with our desired virtual localization
formula,  we see that it suffices to establish:
\begin{equation}
\label{point}
\frac{\Xvir \cdot [Y_i]} { e(TY^m)} = 
\frac{ \Xivir \cap e(E_1^m)} { e(E_0^m)}
\end{equation} 
in $A^{\com^*}(X_i) _t$.
The refined
intersection of a basic  linear
equivalence due to Vistoli [V] with the zero section of a bundle
will yield equation (\ref{point}). The method follows similar
arguments in [BF].

We first review Vistoli's rational equivalence.
Consider the following Cartesian diagram:
\begin{equation}
\label{carty}
\begin{CD}
 \iota^* C_{X/Y}  @>>> C_{X/Y} \\
   @VVV  @ VVV\\
 X_i @>>> X \\
 @VVV @VVV \\
 Y_i @>{\iota}>> Y 
\end{CD}
\end{equation}
The cone 
$C_{X_i/Y_i}$ naturally embeds in
$\iota^* C_{X/Y}$.
Vistoli [V] has constructed a  
rational equivalence in $N_{Y_i/Y} \times \iota^* C_{X/Y}$ which
implies 
\begin{equation}
\label{dwdw}
\iota^![C_{X/Y}]= [C_{X_i/Y_i}]
\end{equation}
in $A_*(\iota^* C_{X/Y})$ (see [BF]). Applying 
Vistoli's equivalence to the $\com^*$-homotopy
quotients yields equation (\ref{dwdw}) in
$A^{\com^*}_*(\iota^* C_{X/Y})$.
We will consider the pull-back of this relation to  
$\iota^*D=  \iota^*C_{X/Y} \times E_0$:
\begin{equation}
\label{vist}
\iota^![D] = 
[D_i \times E_0^m]
\end{equation}
in $A_*^{\com^*}(\iota^* D)$.

Consider the pull-back of
the exact sequence of abelian cones (\ref{fff}) to $X_i$:
$$0 \rarr TY \rarr \iota^*C(I/I^2) \times E_0 \rarr \iota^* C(Q) \rarr 0.$$
There is an inclusion $\iota^* C(Q) \subset E_1$.
The natural inclusion $\iota^* D \subset \iota^* C(I/I^2) \times
E_0$ is $TY$-invariant. 
Hence, the quotient cones 
$$\iota^*D/ TY_i \rarr \iota^*D/ TY \subset \iota^* C(Q)$$
exist.
We obtain a three level
Cartesian diagram:
\begin{equation}
\label{cartone}
\begin{CD}
TY @>>>  \iota^*D \\
@VVV @VVV \\
TY/TY_i @>>> \iota^*D/TY_i \\
@VVV @VVV \\
X_i @>{0_{E_1}}>>  E_1 
\end{CD}
\end{equation}
Note that $TY/TY_i= TY^m$.

We now start the derivation of equation (\ref{point}).
The first steps are:
\begin{eqnarray*}
\Xvir \cdot [Y_i] & = & \iota^! s_{TY}^* 0^!_{E_1} [D] \\
& = & s_{TY}^* 0^!_{E_1} \iota^! [D] \\
& = & s_{TY}^* 0^!_{E_1} [D_i \times E_0^m]
\end{eqnarray*}
in $A_*^{\com^*}(X_i)$.
The first equality is by the definition of $\Xvir$. The second is
obtained from the
commutativity of the intersection product. The third follows from
equation (\ref{vist}).

The $TY_i$-action on $\iota^* D$ leaves the cycle
$D_i \times E_0^m$ invariant (since $TY_i$ acts
naturally on $D_i$ and trivially on $E_0^m$).
By definition, $$D_i/TY_i= \Dvir_i.$$ The class
$[D_i\times E_0^m]\in A_*^{\com^*}(\iota^* D)$ is thus 
the pull-back of
$[\Dvir_i \times E_0^m] \in A_*^{\com^*}(\iota^*D/TY_i)$.
Hence,
$$s_{TY}^* 0^!_{E_1} [D_i \times E_0^m] =
  s_{TY^m}^* 0^!_{E_1} [\Dvir_i \times E_0^m].$$

The scheme-theoretic intersection $0^{-1}_{E_1} (\Dvir_i \times E_0^m)$
certainly lies in $TY^m$. 
The map $$\Dvir_i \times E_0^m \rarr E_1$$ is the product of
the inclusion $\Dvir_i \subset E_1^f$ and
the natural map from the obstruction theory $E_0^m \rarr E_1^m$.
We thus observe $0^{-1}_{E_1}(\Dvir_i \times
E_0^m)$ also lies in $E_0^m$. We conclude the existence
of the following diagram:
\begin{equation}
\label{ccccc}
\begin{CD}
0^{-1}_{E_1}(\Dvir_i \times E_0^m) @>>>  E_0^m \\
@VVV @VVV \\
TY^m @>>>  X_i \\ 
\end{CD}
\end{equation}
To proceed, we need a relation among Gysin maps.

\begin{lm} 
\label{tww}
Let $B_0$ and $B_1$ be $\com^*$-equivariant
bundles on $X_i$.
Let $Z$ be a scheme equipped with two equivariant inclusions 
$j_0$, $j_1$
over $X_i$:
\begin{equation}
\begin{CD}
Z @>>>  B_1 \\
@VVV @VVV \\
B_0 @>>>  X_i \\ 
\end{CD}
\end{equation} Let $\zeta \in A^{\com^*}_*(Z)$.
Then,
$$ s^*_{B_0} j_{0*}(\zeta) \cap e(B_1)
= s^*_{B_1} j_{1*}(\zeta) \cap e(B_0) \ \ \in A^{\com^*}_* (X_i).$$ 
\end{lm}
\bpf
Consider the family of inclusions  $j_t:Z \hookrightarrow B_0 \times B_1$
defined for $t \in \com$ by:
$$j_t= (1-t)\cdot j_0 + t\cdot j_1.$$
The existence of this family implies:
$$s_{B_0 \times B_1}^*j_{0*} (\zeta)= s_{B_0 \times B_1}^*j_{1*}(\zeta).$$
This yields the Lemma by the excess intersection formula.
\epf

Applying Lemma \ref{tww} to diagram (\ref{ccccc}) 
and the class $\zeta= 0^!_{E_1}[\Dvir_i \times E_0^m]$
yields:
\begin{equation}
\label{kkkk}
\Xvir \cdot [Y_i] =
s_{E_0^m}^* \Big( 0^!_{E_1} [\Dvir_i \times E_0^m] \Big) \cdot
\frac{e(TY^m)}{e(E_0^m)}.
\end{equation}

The class $0^!_{E_1} [\Dvir_i \times E_0^m]$ is now considered
to lie in $A_*^{\com^*}(E_0^m)$. As this class does not
depend on the bundle map 
\begin{equation}
\label{aaaa}
E_0^m \rarr E_1^m,
\end{equation}
we are free to assume (\ref{aaaa}) is trivial.
Then, the equality
\begin{equation}
\label{kdf}
s_{E_0^m}^* \Big( 0^!_{E_1} [\Dvir_i \times E_0^m] \Big) =
\Xivir \cap e(E_1^m)
\end{equation}
follows easily from the definition of $\Xivir$ and the
excess intersection formula. Equation (\ref{point})
is a consequence of (\ref{kkkk}) and (\ref{kdf}). The proof of
the virtual localization formula is complete.

\section{\bf{The formula for $\proj^r$}}
\label{projj}
We can use the virtual localization formula (\ref{exloc}) to derive
an expression for the higher genus Gromov-Witten invariants of projective
space analogous to the one given for genus 0 invariants in [K].
The additional arguments needed to justify formula (\ref{exloc}) in the
category of Deligne-Mumford stacks for the moduli space
of maps are given in the Appendices.

We first establish our conventions about the torus
action on projective space.  
Let $V=\com^{r+1}$. Let $p_i\in \proj(V)$ be the points
determined by the basis vectors.
Let $\com^*$ act on $V$ with generic weights $-\lambda_0, \ldots
,-\lambda_r$.  Then, the induced action on the tangent space to $\proj(V)$
at $p_i$ has weights $\lambda_i - \lambda_j$ for $j \neq i$.

Let $\T$ be the full diagonal torus acting on $\proj^r$.
Following [K], we can identify the components of the fixed point
locus of the $\T$-action on $\Mgnp$ with certain types of
marked graphs.  Let $f:C \rarr \proj^r$ be a $\T$-fixed stable map.
The image of $C$ is a $\T$-invariant curve in $\proj(V)$, and  the images
of all marked points, nodes, contracted components, and ramification
points are $\T$-fixed points.   
The  $\T$-fixed points on $\proj^r$ are 
$p_0 , \ldots, p_r$, and the only
invariant curves are the lines joining the points
$p_i$.  It follows that
each non-contracted component of $C$ must map onto one of these lines,
and be ramified only over the two fixed points.  This forces such
a component to be rational, and the map restricted to this 
component is completely determined by its degree.

We are  led to identify the components of the fixed locus with
marked graphs.  To an invariant stable map $f$, we associate a marked 
graph $\Gamma$
as follows.
$\Gamma$ has one edge for each non-contracted component. The  edge $e$ is
marked with the 
degree $d_e$ of the map from that component to its image line.
$\Gamma$ has one vertex for each connected component of $f^{-1}(\{p_0,
\ldots, p_r\})$.
Define the labeling map $$i: \text{Vertices} \rarr \{0, \ldots, r\}$$
by $f(v)= p_{i(v)}$.
The vertices have an additional labeling $g(v)$ by the
arithmetic genus of the associated component.   
(Note the component may be a single point,
in which case its genus is 0.)  Finally, $\Gamma$ has $n$ numbered
legs coming from the $n$ marked points.  These legs are attached to the 
appropriate vertex.  An edge is incident to a vertex if the two
associated subschemes of $C$ are incident.  

The set of all invariant
stable maps whose associated graph is $\Gamma$ is naturally identified
with a finite quotient of a product of moduli spaces of 
pointed curves.  
Define:  $$\M_\Gamma = \prod_{{\rm vertices}} \M_{g(v),\val (v)}.$$
$\M_{0,1}$ and
$\M_{0,2}$ are interpreted as points in this product.
Over the Deligne-Mumford stack $\M_\Gamma$, there is
a canonical family $$\pi: \mathcal{C} \rarr \M_\Gamma$$
of $\com^*$-fixed stable maps
to $\proj^r$ yielding a morphism
$$\gamma: \M_\Gamma \rarr \Mgnp.$$
There is a natural automorphism group $\A$ acting
$\pi$-equivariantly on $\mathcal{C}$ and $\M_\Gamma$.
$\A$ is filtered by an exact sequence of groups:
$$ 1 \rarr \prod_{\rm edges} {\Z}/{d_e} \rarr
\A \rarr \Aut(\Gamma) \rarr 1$$
where $\Aut (\Gamma)$ is the 
automorphism group of $\Gamma$ (as a marked graph).
$\Aut(\Gamma)$ naturally acts on $\prod_{\rm edges} \Z/ d_e$
and $\A$ is the semidirect product. The induced 
map:
$$\gamma/ \A : \M_\Gamma / \A \rarr \Mgnp$$
is a closed immersion of Deligne-Mumford stacks.
It should be noted that the subgroup $\prod_{\rm edges} \Z/ d_e$
acts trivially on $\M_\Gamma$ and that $\M_\Gamma / \A$
is nonsingular.
A component of the $\com^*$-fixed stack of $\Mgnp$
is supported on $\M_\Gamma/ \A$. 
The fixed stack will be shown to be nonsingular by analysis
of the $\com^*$-fixed perfect obstruction theory which yields
the Zariski tangent space.
This nonsingularity
is surprising since the moduli stack $\Mgnp$ is
singular.
A generic $\com^*\subset \T$ will have the same fixed point 
loci in $\Mgnp$.
Via this fixed point identification, the virtual
localization formula will relate the Gromov-Witten
invariants of $\proj^r$ to 
integrals
over moduli spaces of pointed curves.

Following [K], we define a flag $F$ of the graph $\Gamma$ to be an 
incident edge-vertex pair $(e,v)$. Define $i(F)=i(v)$. The edge
$e$ is incident to one other vertex $v'$. Define $j(F)=i(v')$. 
Define: $$\omega_F=\frac{\lambda_{i(F)}-\lambda_{j(F)} }{d_e}.$$ 
This
is the
weight of the induced action of $\com^*$ on the tangent space to the rational
component $C_e$
of $C$ corresponding to $F$ at its preimage over
$p_{i(F)}$.  This fact follows from the corresponding calculation on the weight
of the action on the tangent space to the image line, with a factor of
$\frac{1}{d_e}$ coming from the $d_e$-fold ramification of the map at the 
fixed point.

We describe the obstruction theory of $\Mgnp$ restricted
to $\M_\Gamma/ \A$.  
Define sheaves $\mathcal{T}^1$ and $\mathcal{T}^2$ on
$\M_\Gamma/ \A$
via
the cohomology
of the restriction of  the canonical (dual) perfect 
obstruction theory $E_\bullet$ on $\Mgnp$:
\begin{equation}
\label{eggg}
 0 \rarr \mathcal{T}^1 \rarr E_{0,\Gamma}  \rarr E_{1, \Gamma}
 \rarr \mathcal{T}^2 \rarr 0.
\end{equation}
There is a tangent-obstruction 
exact sequence of sheaves on the substack $\M_\Gamma/\A$:
\begin{equation}
\label{tanob}
0 \rarr \Ext^0(\Omega_C(D), \oh_C) 
\rarr H^0(C,f^*TX) \rarr \mathcal{T}^1 \rarr 
\end{equation}
$$ \rarr \Ext^1(\Omega_C(D), \oh_C) \rarr H^1(C,f^*TX) 
\rarr \mathcal{T}^2 \rarr 0.$$
The marked point divisor on $C$ is denoted by $D$.
The 4 terms other than the sheaves $\mathcal{T}^i$ are vector
bundles and are labeled by their fibers. This sequence can be
viewed as
filtering the
deformations of the maps by those 
which preserve the domain curves. It arises via the 
pull-back to $\M_\Gamma/\A$ of a distinguished triangle of complexes
on $\Mgnp$ (see Appendix B).
These results may be found in [LT], [R2], [B].

In the remainder of this section, the fixed and moving
parts of the 4 bundles in the tangent-obstruction complex
are explicitly identified following [K].
It is simpler to carry out the bundle analysis on the
prequotient $\overline{M}_\Gamma$ to avoid monodromy
in the nodes. In fact, the final integrals over the
fixed locus will be evaluated on $\M_\Gamma$
and corrected by the order  of $\A$.

It will be seen that there are exactly 3
fixed pieces in the 4 bundles. 
They occur in the $1^{\rm st}$, $2^{\rm nd}$, and
$4^{\rm th}$ terms of the complex. The fixed piece in the $1^{\rm st}$ term
maps isomorphically to the fixed piece of the $2^{\rm nd}$.
$\mathcal{T}^{1,f}$ is thus isomorphic to
the fixed piece in the $4^{\rm th}$ term. The latter  is 
canonically the tangent bundle to $\M_\Gamma$.
Also, $\mathcal{T}^{2,f}=0$. We can conclude that the
fixed stack is nonsingular and equal to $\M_\Gamma/\A$. 
The two exact sequences 
(\ref{eggg}) and (\ref{tanob}) imply:
$$e(\Nvir)= \frac{e(B_2^{m}) e(B_4^{m})}{e(B_1^m)e(B_5^m)}$$
where, for example, $B_2^m$ denotes the moving part
of the $2^{\rm nd}$ term in (\ref{tanob}).

We first calculate the contribution coming from the
bundle $$\Aut(C)=\Ext(\Omega_C(D), \oh_C)$$ parameterizing infinitesimal
automorphisms of the pointed domain.
For each non-contracted component of $C$,
there is  a weight zero piece coming
from the infinitesimal automorphism of that component
fixing the two special points.  This term will cancel with a similar term
in $H^0(f^*T\proj^r)$.
Also, since there is no moving part, $e(B_1^{m})=1$.   
If it is the case that the special points are not
marked or nodes, that is the associated vertex of the graph has genus 0 and
 valence
one, there would be an extra automorphism with nontrivial weight.  We 
will leave this case and the case of a genus 0 valence 2
vertex to the reader. 
No extra trivial weight pieces arise in these two
cases.  
As in [K], the (integrated) final formulas for the genus 0 vertex
contributions 
will still be correct.

Next, we consider the bundle $\Def(C)=\Ext^1(\Omega_C(D), \oh_C)$ 
parameterizing
deformations of the pointed domain.
A deformation of the contracted components (as marked curves) is
a weight zero deformation of the map which yields the
tangent space of $\M_\Gamma/ \A$
as a summand 
in the weight zero piece
of $\Def(C)$.  
The other deformations of $C$ come from smoothing nodes of $C$ which join
contracted components to non-contracted components.  This space splits
into a product of spaces corresponding to deformations which smooth
each node individually.  The one dimensional space associated to each
node is identified as a bundle with the tensor product of the tangent
spaces of the two components at the node.  We  see
that the tangent space to the non-contracted curve forms a trivial bundle
with weight $\omega_F$ while the tangent space to the contracted curve
varies but has trivial weight.  Let  $e_F$  denote the line
bundle on $\M_\Gamma$ whose fiber over a point is the cotangent space
to the component associated to $F$ at the corresponding node.  Therefore,
$$e(B_4^{m})= \prod_{\rm flags}(\omega_F - e_F).$$

To compute the contribution coming from $H^\bullet(f^*T\proj^r)$, we
consider the normalization sequence resolving all of the nodes
of $C$ which are forced by the graph type $\Gamma$.  
$$0 \rarr \oh_C \rarr \bigoplus_{\rm vertices} \oh_{C_v} \oplus 
\bigoplus_{\rm edges} \oh_{C_e}
\rarr \bigoplus_{\rm flags} \oh_{x_F} \rarr 0$$
Twisting by $f^*(T\proj^r)$ and taking cohomology yields:
$$0 \rarr H^0(f^*T\proj^r) \rarr  
\bigoplus_{\rm vertices} H^0(C_v, f^*T\proj^r) 
\oplus \bigoplus_{\rm edges}H^0(C_e, f^*T\proj^r) \rarr$$
$$\rarr \bigoplus_{\rm flags}
T_{p_{i(F)}}\proj^r \rarr H^1(f^*T\proj^r) \rarr
 \bigoplus_{\rm vertices} H^1(C_v, f^*T\proj^r) 
\rarr 0$$
where we have used the fact that there will be no higher cohomology on
the non-contracted components since they are rational.
Also note that $H^0(C_v, f^*(T\proj^r) = T_{p_{i(v)}}\proj^r$
since $C_v$ is connected and $f$ is constant on it.
Thus, we obtain:
$$
H^0 - H^1 =  \begin{array}{cccc} + & {\displaystyle \bigoplus_{\rm vertices} 
T_{p_{i(v)}} 
\proj^r}& +&{\displaystyle \bigoplus_{\rm edges} H^0(C_e, f^*T\proj^r)}   \\
- & {\displaystyle \bigoplus_{\rm flags} T_{p_{i(F)}}\proj^r}& -& 
{\displaystyle \bigoplus_{\rm vertices} H^1(C_v, f^*T\proj^r)}
\end{array}
$$
As non-contracted components are rigid, we see that $H^0(C_e,f^*T\proj^r)$
is trivial as a bundle, but we need to determine its weights.  We do this
via the Euler sequence.
On $\proj^r$ we have:
$$0 \rarr \oh \rarr \oh(1) \otimes V \rarr T\proj^r \rarr 0. $$
Pulling back to $C_e$ and taking cohomology gives us:
$$ 0 \rarr \C \rarr H^0(\oh(d_e))\otimes V \rarr H^0(f^*T\proj^r) \rarr 0 $$
Here the weight on $\C$ is trivial, and the weights on $H^0(\oh(d_e))$
are given by $\frac{a}{d_e}\lambda_i + \frac{b}{d_e}\lambda_j$ for $a+b=d_e$.
The weights on $V$ are $-\lambda_0, \ldots ,-\lambda_r$.  So the
weights of the middle term are just the pairwise sums of these, 
$\frac{a}{d_e}\lambda_i + \frac{b}{d_e}\lambda_j -\lambda_k$.  There are
exactly 2 zero weight terms here coming from $a=0, k=j$ and $b=0, k=i$.  
These cancel the zero weight term from the $\C$ on the left, and the zero
weight term occurring in $\Aut(C)$.  Breaking up the remaining terms into
two groups corresponding to $k=i,j$ and $k \neq i,j$, we obtain
the contribution
of $H^0(C_e,f^*T\proj^r)$ to the Euler class ratio $e(B_2^{m})/e(B_5^{m})$:
$$(-1)^{d_e} \frac{{d_e!}^2 }{ d_e^{2d_e}}
(\lambda_i - \lambda_j)^{2d_e}
\cdot \prod_{\stackrel{a+b=d_e}{ k\neq i,j}} (\frac{a} {d_e}\lambda_i + 
\frac{b}{d_e}\lambda_j -\lambda_k).$$

Finally, we evaluate the contribution of $H^1(C_v, f^* T\proj^r)$.  This is 
simply $H^1(C_v, \oh_{C_v})\otimes T_{p_{i(v)}}\proj^r$.  
As a bundle, $H^1(C_v, \oh_{C_v})$ is the dual of the Hodge bundle 
$E=\pi_* \omega$ on the moduli space $\M_{g(v), {\rm val}(v)}$.
The bundle $H^1(C_v, \oh_{C_v})\otimes T_{p_{i(v)}}\proj^r$
splits into $r$ copies of $E^\vee$ twisted respectively
by the $r$ weights $\lambda_i -
\lambda_j$ for $j \neq i$.  Taking the equivariant top Chern class
of this bundle yields:
$$\prod_{j \neq i} c_{(\lambda_i - \lambda_j)^{-1}}(E^\vee)\cdot 
(\lambda_i - \lambda_j)^{g(v)}$$
where for a bundle $Q$ of rank $q$:
$$c_t(Q)= 1+ t c_1(Q) + \ldots t^q c_{q}(Q).$$

We arrive at the following form
of the inverse Euler class of the virtual normal bundle
to the fixed point locus corresponding to the graph $\Gamma$. 

\begin{eqnarray*}
& &\prod_{\rm flags} \frac{1}{\omega_F - e_F} 
\prod_{j \neq i(F)} (\lambda_{i(F)} -\lambda_j) \\
 \frac{1}{e(\Nvir)}& =&\prod_{\rm vertices} \prod_{j \neq i(v)}
c_{(\lambda_{i(v)} - \lambda_j)^{-1}}(E^\vee) \cdot (\lambda_{i(v)}
-\lambda_j)^{g(v)-1} \\
& &\prod_{\rm edges} \frac{(-1)^{d_e} d_e^{2d_e}} { (d_e!)^2 
(\lambda_i -\lambda_j)^{2d_e}}
\prod_{\stackrel{a+b=d_e} {k\neq i,j}} 
\frac{1}{ \frac{a}{d_e}\lambda_i + \frac{b}{d_e}\lambda_j
-\lambda_k}
\end{eqnarray*}

In addition, the virtual
fundamental class of the fixed locus must be identified.
We have already seen $\mathcal{T}^{1,f}$ is
the tangent bundle of $\M_\Gamma$.
and $\mathcal{T}^{2,f}=0$. It then follows from
(\ref{eggg}) that
the $\com^*$-fixed (dual) perfect obstruction
theory is equivalent  on the fixed stack
to the trivial perfect obstruction theory. The
virtual fundamental class of the fixed stack is simply
the ordinary fundamental class.

The above expression for $\frac{1}{e(\Nvir)}$
can be used in the virtual 
localization formula to deduce formulas expressing
Gromov-Witten invariants of projective space in terms
of integrals on moduli spaces of pointed curves.  The numerator
terms, coming from the cohomology classes of $\proj^r$
are identical in this
higher genus case to the terms appearing in [K].  In particular,
they contribute only additional weights, and no cohomological
terms.

Let $[n]=\{1,\ldots,n\}$ be the marking set of an $n$-pointed
graph $\Gamma$.
Let $i:[n] \rarr \{0, \ldots, r\}$ be defined
by $f(m)= p_{i(m)}$. The final expression
for the Gromov-Witten invariants of $\proj^r$ is: 
$$ I_{g, d}^{\proj^r}(H^{l_1}, \ldots, H^{l_n})=
\sum_{\Gamma} 
\frac{1}{|\A_\Gamma|}
\int_{\M_\Gamma} \frac{\prod_{[n]} \lambda_{i(m)}^{l_m}}
{e(\Nvir_\Gamma)}.$$
The sum is over all graphs $\Gamma$ indexing fixed
loci of $\Mgnp$.
To evaluate the integral, one expands the terms of the form
$\frac{1}{\omega - e}$ as formal power series, and then integrates
all terms of the appropriate degree.

Each integral that is encountered will naturally split as a product
of integrals over the different moduli spaces of pointed curves.
We remark that the integrals over genus 0 spaces are identical
to the ones which are dealt with in [K].  In particular, while the
formula given above is incorrect for graphs with vertices of genus
0 and valence 1 or 2, the formulas obtained in [K] after integrating
over $\M_{0,n}$ hold for these degenerate cases as well.

In higher genera, we know of no closed formulas for the integrals
which occur in these calculations.  However, C. Faber has constructed
an algorithm in [Fa] which determines all such integrals.  Thus, this
formula, together with Faber's algorithm, 
gives a method in principle to determine arbitrary Gromov-Witten
invariants of projective space.

\section{\bf{Multiple cover calculations}}
\label{ellip}
Let $C\subset X$ be
a nonsingular rational curve with balanced normal
bundle $N\eqq\oh(-1) \oplus \oh(-1)$
in a nonsingular Calabi-Yau 3-fold $X$.
Let $[C]\in H_2(X, \mathbb{Z})$ be the
homology class of $C$. 
The space of stable elliptic maps to $X$ 
representing the curve class $d[C]$ contains
a component $Y_d$ 
consisting of maps which factor through
a $d$-fold cover of $C$. $Y_d$ is naturally
isomorphic to 
$\barr{M}_{1,0}(C,d)$, the space
of unpointed, genus 1 stable maps. 
The contribution
of $Y_d$ to the elliptic Gromov-Witten invariant
$I_{1,d[C]}^X$ has been computed in physics [BCOV].
The answer obtained is $\frac{1}{12d}$ (accounting
for the differing treatment of the elliptic involution).

Mathematically, the excess contribution of $Y_d$
is expressed as an integral over $\barr{M}_{1,0}(C,d)$.
The integral is computed here for all $d$ via localization. Localization
reduces the contribution to a graph sum which can
be explicitly evaluated by Manin's trick [M] and
a formula for intersections of cotangent
lines on $\overline{M}_{1,n}$.

Let $\pi: U \rarr \barr{M}_{1,0}(C,d)$ be the
universal family over the moduli space. Let
$\mu: U \rarr C$
be the universal evaluation map.
The expected dimension
of $\barr{M}_{1,0}(C,d)$ is $2d$. 
By the cohomology and base change theorems,
$R^1 \pi_* \mu^* N$ is a vector bundle of rank
$2d$ on $\barr{M}_{1,0}(\proj^1,d)$. 
The contribution of $Y_d$ to the elliptic Gromov-Witten
invariant of curve class $d[C]$ is:
\begin{equation}
\label{exx}
\int _{[ \barr{M}_{1,0}(C,d)]^{\it vir}} c_{2d} (R^1 \pi_* \mu^* N).
\end{equation}
Natural lifts of $\com^*$-actions on $C$ to
$\barr{M}_{1,0}(C,d)$,
$N$, and $R^1 \pi_* \mu^* N$ exist. The localization formula
can therefore be applied to compute (\ref{exx}).
The answer obtained agrees with the physics calculation.
\begin{pr}
\label{mule}
$$\int _{[ \barr{M}_{1,0}(C,d)]^{\it vir}} c_{2d} (R^1 \pi_* \mu^* N) =
\frac{1}{12 d}.$$
\end{pr}

Let $V\eqq \com^2$. Let $C= \proj(V)$. 
Let $\com^*$ act by weights
$0$ and $-1$ on $V$. Let $x_0$ and $x_{-1}$ be the respective 
fixed points in $C$. 
The $\com^*$-action lifts naturally to the tautological
line $\oh(-1)$ and thus to $N$. 
Consider the graph sum obtained by the localization
formula for the integral (\ref{exx}).
The 0 weight leads
to a drastic collapse of the sum. This was 
observed by Manin in [M] for an analogous excess
integral over a space of genus 0 maps.
In fact, the only graphs which contribute
are comb graphs where the
backbone is an elliptic curve contracted over
$x_{-1}$ and the teeth are rational curves 
multiple covering $\proj(V)$.
The degree $d$ is distributed over the
teeth by $\sum_{1}^{k} m_i =d$.
The denominator terms in the localization formula are determined
by the results of Section \ref{projj}. The numerator is given by the bundle
$R^1\pi_* \mu^*N$ which is decomposed on each fixed
point locus via the natural normalization sequence.
The formula 
\begin{equation}
\label{mast}
\sum_{m\vdash d} \frac{(-1)^{d-L(m)}}{\Aut(m) \ \Pi_{1}^{L(m)} m_i}
\int_{\overline{M}_{1,L(m)}} \frac{1+\lambda}{\Pi_{1}^{L(m)}
(1-m_ie_i)}
\end{equation}
is obtained for the degree $d$ contribution. 
The sum is over all positive partitions:
$$m=(m_1, \ldots, m_k), \ \ m_i>0, \ \ \sum_{1}^{k} m_i =d.$$  
$L(m)$ denotes
the length of $m$. $\Aut(m)$ is
the order of the stabilizer of the symmetric group
$S_k$-action on the string $(m_1, \ldots, m_k)$.
The class $\lambda$ in the numerator is the first
Chern class of the Hodge bundle on $\overline{M}_{1,n}$.
As before, $e_i$ is the $i^{th}$ cotangent line bundle
on $\overline{M}_{1,n}$.

The integral (\ref{mast}) is calculated in two parts to
prove Proposition \ref{mule}.
\begin{lm}
\label{AAA}
$$
\sum_{m\vdash d} \frac{(-1)^{d-L(m)}}{\Aut(m) \ \Pi_{1}^{L(m)} m_i}
\int_{\overline{M}_{1,L(m)}} \frac{\lambda}{\Pi_{1}^{L(m)}
(1-m_ie_i)}= \frac{d}{24}$$
\end{lm}
\begin{lm}
\label{BBB}
$$
\sum_{m\vdash d} \frac{(-1)^{d-L(m)}}{\Aut(m) \ \Pi_{1}^{L(m)} m_i}
\int_{\overline{M}_{1,L(m)}} \frac{1}{\Pi_{1}^{L(m)}
(1-m_ie_i)}= \frac{d}{24}$$
\end{lm}

We start with Lemma \ref{AAA}.
The first step is to use the boundary expression for $\lambda$ to
reduce to an integral over genus 0 pointed moduli
spaces. On $\overline{M}_{1,1}$, the equation:
\begin{equation}
\label{lamb}
\lambda = \frac{\Delta_0}{12}
\end{equation}
holds where $\Delta_0$ is the irreducible boundary divisor.
Since $\lambda$ on $\overline{M}_{1,n}$ is a pull-back from
a one pointed space,
(\ref{lamb}) is valid on $\overline{M}_{1,n}$.
Using the standard identification of $\Delta_0$ with the
$\mathbb{Z}/2\mathbb{Z}$-quotient of $\overline{M}_{0,n+2}$, the
equality:
$$\int_{\overline{M}_{1, L(m)}} \frac{\lambda}{\Pi_{1}^{L(m)}
(1-m_ie_i)} = \frac{1}{24} \int_{\overline{M}_{0, L(m)+2}}
\frac{1}{\Pi_{1}^{L(m)}
(1-m_ie_i)}$$
is obtained.
Next, using the well-known formula for 
intersection numbers on the genus 0 spaces, we see:
$$\int_{\overline{M}_{0, L(m)+2}}
\frac{1}{\Pi_{1}^{L(m)}
(1-m_ie_i)} = (\sum_{1}^{L(m)} m_i)^{L(m)-1}= d^{L(m)-1}.$$
After substituting these equalities, the sum of
Lemma \ref{AAA} is transformed to:
\begin{equation}
\label{manin}
\frac{(-1)^d}{24d} \sum_{m\vdash d} \frac{(-1)^{-L(m)}}{\Aut(m)\  
\Pi_1^{L(m)}m_i}
{d}^{L(m)}.
\end{equation}
The summation term 
in (\ref{manin}) was encountered by Manin in [M]. It evaluates
explicitly to $(-1)^d$ via a  generating function
argument (see [M] p.416). The value of (\ref{manin}) is thus
$\frac{1}{24d}$.
Lemma \ref{AAA} is established.

We now prove Lemma \ref{BBB}. A generating function
approach is taken. For $d\geq 1$, let
$$g_d= 
\sum_{m\vdash d} \frac{(-1)^{d-L(m)}}{\Aut(m) \ \Pi_{1}^{L(m)} m_i}
\int_{\overline{M}_{1,L(m)}} \frac{1}{\Pi_{1}^{L(m)}
(1-m_ie_i)}.$$
Define $\gamma(t)$ by:
$$\gamma(t)= \sum _{\alpha \geq 1} (-1)^{\alpha} g_\alpha t^\alpha.$$
An important observation is that $\gamma(t)$ can be
rewritten in the
following form:
\begin{equation}
\label{oboe}
\gamma(t)= < 
\exp \ (-\sum_{\alpha\geq 1} \sum_{i\geq 0}\alpha^{i-1} t^\alpha
\sigma_i) >_1.
\end{equation}
Here, Witten's notation,
\begin{equation}
\label{jenny}
<\sigma_{0}^{r_0} \sigma_{1}^{r_1} \cdots \sigma_{k}^{r_k}>_1, 
\end{equation}
is used to denote the integral: 
$$\int_{\overline{M}_{1,r}}
\underbrace {e_{r_0+1} \ldots e_{r_0+r_1}}_{r_1} \cdot
\underbrace {e^2_{r_0+r_1+1} \ldots
e^2_{r_0+r_1+r_2}}_{r_2} \ldots  
\underbrace
{e^k_{r-r_k+1} \ldots
e^k_{r}} _{r_k}$$
where $r=\sum_1^k r_i$.
Equality (\ref{oboe}) is a simply a rewriting of terms.
The genus 1 integrals (\ref{jenny}) are determined
from genus 0 integrals
by a beautiful formula in the formal variables $\{z_i\} _{i \geq 0}$:
\begin{equation}
\label{dike}
< \exp \sum_{i\geq 0} z_i \sigma_i >_1 =
\frac{1}{24} \log < \sigma_0^3 \exp
\sum_{i\geq 0} z_i \sigma_i >_0 .
\end{equation}
Formula (\ref{dike}) can be found, for example, in [D].
Let $z_i=-\sum_{\alpha\geq 1} \alpha^{i-1}t^\alpha$.
Using (\ref{oboe}) and (\ref{dike}), $\gamma(t)$ may
be expressed as:
\begin{equation}
\label{gam}
\gamma(t)=\frac{1}{24} \log  < \sigma_0^3 \exp
(-\sum_{\alpha\geq 1} \sum_{i\geq 0}\alpha^{i-1} t^\alpha
\sigma_i) >_0
\end{equation}

Equation (\ref{gam}) will be used to determine $\gamma(t)$.
First, 
define another generating function $\psi(t)$ by:
$$\psi(t)= 1 + \sum_{\beta} s_\beta t^\beta$$
where the coefficients $s_\beta$ are:
\begin{equation}
\label{quail}
s_\beta= \sum_{m\vdash \beta} \frac{ (-1)^{-L(m)}}{\Aut(m) \ 
\Pi_{1}^{L(m)} m_i} \int _{\overline{M}_{0, L(m)+3}} \frac{1}
{\Pi_{1}^{L(m)} (1-m_ie_i)}.
\end{equation}
As before,  the equality:
$$\psi(t)= < \sigma_0^3 
\exp(-\sum_{\alpha\geq 1} \sum_{i\geq 0}\alpha^{i-1} t^\alpha
\sigma_i) >_0$$
is a rewriting of terms.
However, the expression (\ref{quail}) may be explicitly
evaluated by the genus 0 intersection formulas and
Manin's summation argument to yield:
$$ s_\beta= (-1)^\beta.$$
Hence, $\psi(t)$ is simply $1/(1+t)$, and
$$\gamma(t)= -\frac{\log(1+t)}{24}= \frac{1}{24}(-t + \frac{t^2}{2}-
\frac{t^3}{3} + \ldots).$$
Thus, $g_d= \frac{1}{24d}$. Lemma \ref{BBB} is proven.
Proposition \ref{mule} follows from (\ref{mast}) and the
two Lemmas.

Localization may be applied to the analogous excess 
integrals for arbitrary genus $g$. The resulting
formula is:

\begin{equation}
\label{mast2}
\int _{[ \barr{M}_{g,0}(\proj^1,d)]^{\it vir}} 
c_{\rm {top}} (R^1 \pi_* \mu^* N) =
\end{equation}
\vspace{+15pt}
$$
\sum_{m\vdash d} \frac{(-1)^{d-L(m)}}{\Aut(m) \ \Pi_{1}^{L(m)} m_i}
\int_{\overline{M}_{g,L(m)}} \frac{1+c_1(E)+\ldots + c_g(E)}{\Pi_{1}^{L(m)}
(1-m_ie_i)}$$
\vspace{+15pt}

\noindent where $E$ is the Hodge bundle.
For $g\geq 2$, we have conjectured 
 with C. Faber the above integral sum is equal to:
\begin{equation*}
\frac{|B_{2g}| \cdot d^{2g-3}}{2g\cdot (2g-2)!}
= \frac{|\chi(M_g)|\cdot d^{2g-3}}{(2g-3)!}.
\end{equation*}
This equality has been verified in case $g+d \leq 7$.

\appendix

\section{\bf{Global nonsingular embeddings}}
Let $V$ be an nonsingular projective
algebraic variety with a $\com^*$-action. After a group cover
$\com^*\rarr \com^*$, a $\com^*$-equivariant polarization $\mathcal{L}$
on $V$ may be found (see [MFK]).  
$\overline{M}_{g,n}(V,\beta)$ is a $\com^*$-equivariant closed
substack of $\overline{M}_{g,n}(\proj^r,\beta)$ via the
equivariant embedding determined by $\mathcal{L}$.
It will be shown that the Deligne-Mumford stack
$\overline{M}_{g,n}(\proj^r,\beta)$ admits a global closed
equivariant embedding
into a nonsingular Deligne-Mumford stack.

In [FP], the moduli space of
maps to $\proj^r$ is expressed as a quotient of a locally
closed scheme $J$ of an
associated (product)
Hilbert scheme $\mathcal{H}$ by a reductive group $\G= \mathbf{PGL}$.
Four properties of this quotient will
be needed here.
\begin{enumerate}
\item[(i)] The stack quotient $[J/\G]$ is the
Deligne-Mumford moduli stack of maps $\overline{M}_{g,n}(\proj^r,\beta)$.
\item[(ii)] $\G$ acts with finite stabilizers on $J$ (in fact, the
    $\G$-action on $J$ is proper).
\item[(iii)] There is a $\com^*\times \G$-action on $J$
             which descends to the given $\com^*$-action
             on $\overline{M}_{g,n}(\proj^r,\beta)$.
\item[(iv)] There is a $\com^*\times \G$-equivariant linearized
            embedding
            of $J\subset \mathcal{H}$ 
             in a nonsingular Grassmannian $\mathbb{G}$.
\end{enumerate}
All these properties are obtained directly from the
construction in [FP].  

The $\G$-equivariant open set $U\subset \mathbb{G}$ on
which the $\G$-action has finite stabilizers contains
$J$ and is $\com^*$-equivariant. 
Note that $\Delta=\barr{J} \setminus J$ is closed 
in $\mathbb{G}$ and is $\com^* \times \G$-equivariant.
After 
discarding $\Delta \cap U$,
it may be assumed that $J$ is closed in $U$. 
Let $Y$ be the nonsingular Deligne-Mumford stack $[U/\G]$. 
The moduli space of maps embeds $\com^*$-equivariantly in $Y$.
It should be noted that while $Y$ is a (quasi-separated)
Deligne-Mumford stack of finite type, $Y$ need not
be separated.

\section{{\bf The obstruction theory on} $\overline{M}_{g,n}(V,\beta)$}
In [B] and [LT], a canonical 
obstruction theory on $\overline{M}_{g,n}(V,\beta)$ is defined
which is locally a two term complex of vector bundles.
To obtain a global two term complex, a polarization is required. 
Since the constructions are $\com^*$-equivariant, 
a $\com^*$-equivariant perfect obstruction theory on
$\overline{M}_{g,n}(V,\beta)$ may be defined using $\mathcal{L}$.

We sketch here the method of [B] to obtain the
equivariant perfect obstruction theory on 
$\overline{M}_{g,n}(V,\beta)$.
The relative deformation problem is considered
for the canonical morphism 
$$\tau: 
\overline{M}_{g,n}(V,\beta) \rarr \mathfrak{M}_{g,n}$$
where $\mathfrak{M}_{g,n}$ is the nonsingular Artin
stack of prestable curves.
The theory of the cotangent complex for Artin stacks
has been developed in [LM-B]. The morphism $\tau$ yields
a distinguished triangle of cotangent
complexes on $\overline{M}_{g,n}(V,\beta)$:
$$\tau^* L_{\mathfrak{M}}^\bullet \rarr L_{\barr{M}}^\bullet
 \rarr L_\tau^\bullet \rarr  
\tau^* L_{\mathfrak{M}} ^\bullet [1].$$
The complex $\tau^* L_{\mathfrak{M}}^\bullet$ has a two term 
bundle resolution
of amplitude $[0,1]$, $A^\bullet=A^0 \rarr A^1$,
in the $\com^*$-equivariant derived category (obtained from
$\mathcal{L}$). 
There is a relative obstruction theory (see [B])
\begin{equation}
\label{relly}
R^\bullet \pi_*(f^*T_V)^{\vee} \rarr L_\tau^\bullet
\end{equation}
with a natural map to $\tau^* L_{\mathfrak{M}}^\bullet [1]$.
Here, $f$ is the map to $V$ from the universal
curve and $\pi$ is the projection from the universal curve
to $\barr{M}_{g,n}(V, \beta)$. Moreover, $R^\bullet\pi_*(f^*T_V)^{\vee}$
has a two term equivariant bundle resolution $B^\bullet=
B^{-1} \rarr B^0$. Representatives of $L_\tau ^\bullet$ and
$L_{\barr{M}}^\bullet$
in the equivariant derived category may be found from
the global nonsingular embeddings constructed above.
The diagram below of distinguished triangles may be
canonically completed in the equivariant derived
category to obtain an equivariant
obstruction theory of amplitude [-1,1]:
$$E^\bullet= B^{-1} \rarr B^0 \oplus A^0 \rarr A^1,$$ 
\begin{equation}
\label{hdm}
\begin{CD}
B^\bullet  @>>> A^\bullet [1] @>>> E^\bullet [1]  @>>> B^\bullet [1]\\
@VVV   @VVV  @VVV @VVV\\
L^\bullet_\tau  @>>>\tau^* L_{\frak{M}}^\bullet [1]
@>>> L_{\barr{M}}^\bullet [1] @>>> L^\bullet_\tau [1].
\end{CD}
\end{equation}
The stability condition implies the cohomology
of $E^\bullet$ at grade 1 vanishes. Hence, $E^\bullet$
may be represented by a two term equivariant complex
$E^{-1} \rarr E^{0}$. 
The morphism constructed in diagram (\ref{hdm}),
$$\phi: E^\bullet \rarr L^\bullet_{\barr{M}}$$
can be seen to be a equivariant perfect obstruction
theory.

Property (iv) of Appendix A implies that the moduli stack 
$\overline{M}_{g,n}(V,\beta)$
has enough $\com^*$-equivariant locally frees. Hence,
representatives  $(E^\bullet,\phi)$ in the equivariant
derived category may be chosen to map to
the two term cutoff $[I/I^2 \rarr \Omega_Y]$ of
$L^\bullet_{\barr{M}}$
determined by the nonsingular embedding (see Section \ref{locfor}).
The distinguished triangle (\ref{hdm}) is used in the
computations of Section \ref{projj}.

\section{\bf{Localization for Deligne-Mumford stacks}}
In this appendix, we extend the virtual localization
formula to the case of a Deligne-Mumford
stack $X$ with a $\C^*$-equivariant perfect obstruction
theory under the additional assumption that
$X$
admits an equivariant global embedding in a nonsingular
Deligne-Mumford stack.
This condition is satisfied for the
moduli space of maps to a nonsingular projective variety $V$
with a $\com^*$-action by the 
existence of a quotient construction reviewed in Appendix A.
We assume Deligne-Mumford stacks are of finite type
(but not necessarily separated). In fact,
for the application to the moduli space of maps, only
quotient stacks $[U/\G]$ (where $U$ is a quasi-projective
scheme and $\G$ acts with finite stabilizers) need be considered.
The $\com^*$-actions on  $[U/\G]$ 
which arise descend from $\com^* \times \G$-actions
on $U$.

First of all, a good theory of rational Chow groups
on Deligne-Mumford stacks has been constructed in
[V]. 
A finite covering result of [LM-B] (Theorem 10.1) is required
for Vistoli's theory to apply to  general Deligne-Mumford stacks. 
Essentially all properties
of Chow groups for schemes hold for Deligne-Mumford
stacks.  In particular, one obtains flat pullbacks
and proper pushforwards.  Flat pullback gives
an isomorphism in Chow groups between any stack and
any vector bundle over the stack.
Refined Gysin maps exist for regular embeddings,
giving rise to an intersection product on
the Chow groups of smooth Deligne-Mumford stacks.  
Finally, Chern and Segre classes
for vector bundles and cones exist and satisfy
the same properties as they do for schemes. 
Also, because these groups are defined
in terms of closed substacks, it is immediate that the Chow groups
are non-zero only in dimensions
between zero and the dimension of the 
stack.  (This last condition would not be possible if one required
a theory with integral
coefficients satisfying the hypothesis that
flat pullback to a vector bundle gave isomorphisms
in Chow groups.)  

We can define a
theory of $\com^*$-equivariant Chow groups on 
Deligne-Mumford stacks by the following constructions.  
One simply defines $A^{\com^*}_* (X)$ to be the
Chow groups of appropriate approximations
to the homotopy quotient as in [T], [EG1].  The only difference
is that the homotopy quotient is now taken to be the
stack quotient, $[X \times E\com^*/ \com^*]$. 
Since this is a free group action,
the quotient is also a Deligne-Mumford stack.
Therefore, $\com^*$-equivariant Chow
groups are well defined for Deligne-Mumford stacks:
$$A^{\com^*}_* (X) = A_*  [X \times E\com^*/ \com^*].$$
The proof in Section \ref{gencase} proves the 
virtual localization formula on Deligne-Mumford
stacks satisfying the embedding assumption provided we
know that the standard localization formula 
holds for nonsingular Deligne-Mumford stacks.

The key step in the proof the localization formula
for nonsingular Deligne-Mumford stacks, 
as in proofs of localization
in other categories, is the following lemma:
\begin{lm}
If $U$ is a Deligne-Mumford stack on which
$\com^*$ acts without fixed points, then
the equivariant Chow groups $A^{\com^*}_*(U)$ vanish
after localization.
\end{lm}

\bpf
We consider the stack quotient, $[U/\C^*].$  Because
the $\com^*$-action has no fixed points on $U$, 
the quotient is again Deligne-Mumford.
Furthermore,
we claim that the Chow groups of $[U/\com^*]$ 
are naturally isomorphic
to the equivariant Chow groups of $U$.  
In the diagram below, both horizontal arrows
are open sets of vector bundles.
\begin{equation*}
\begin{CD}
U \times E\com^*_k  @>>> U \\
@VVV   @VVV \\
 [U \times E\com^*_k/ \com^*]  @>>> [U/\com^*]
\end{CD}
\end{equation*}
\noindent $E\com^*_k$ is an approximation
to $E\com^*$ determined by an open set of a
$\com^*$-representation ([T], [EG1]).
As these approximations to the 
homotopy quotient are realized as 
open sets of vector bundles over $[U/\C^*]$, the
isomorphism in Chow groups follows.
We have already observed that $A_* [U/\C^*]$
has only finitely many graded components.  Hence, $A_* ^{\com^*}(U)$
has only finitely many graded components and thus
is  trivial after localization.
\epf

Let $Y$ be a nonsingular Deligne-Mumford stack with a
$\com^*$-action. The fixed substack $Y^f$ is defined
as the stack theoretic zero locus of the canonical
vector field determined in $TY$ by the
flow. 
It follows easily that $\com^*$ has no fixed points
on $U=Y\setminus Y^f$.
If we consider the pushforward
$\iota_*$
from the equivariant Chow groups of the fixed locus,
$Y^f$ to $Y$, we obtain the exact sequence in equivariant
Chow groups:
$$ A^{\C^*}_*(Y^f) \stackrel{\iota_*}{\rarr} 
A^{\C^*}_*(Y) \rarr A^{\C^*}_*(U) \rarr 0.$$ 
By the Lemma, $\iota_*$ is surjective
after localization.

To prove $\iota_*$ is injective after localization,
the nonsingularity of $Y^f$  will
be used.
Some care must to be taken to establish the nonsingularity
in the category of Deligne-Mumford stacks.
Let $\zeta$ be a $\com^*$-fixed geometric point.
The stack analogue of the local ring is the
strict Henselization  $H_\zeta$ of the
local ring of $\zeta$ determined in an \'etale neighborhood.
$\com^*$ (more precisely
 a finite cover $\com^* \rarr \com^*$ of the
original $\com^*$) acts equivariantly on all
quotients of $H_\zeta / m_\zeta^p$ where $m_\zeta$
is the maximal ideal. The ideal $J_\zeta \subset H_\zeta$
of the $\com^*$-fixed scheme at $\zeta$ may be analyzed
via the $\com^*$-actions on these quotients. 
In fact, R. Thomason in [T] p.456 has proven in this
case $J_\zeta$ is generated via a regular sequence
in regular local ring $H_\zeta$. It follows that the
$\com^*$-fixed stack $Y^f$ is nonsingular. 

The $\com^*$-action on $Y^f$ is equivalent to the
trivial $\com^*$-action after a finite cover $\com^*\rarr \com^*$
of the original $\com^*$. Then,
$\com^*$-equivariant sheaves on $Y^f$
split into fixed and moving parts as before.
The normal bundle to
the fixed locus
is seen to have the standard identification with the
moving part of the restriction of the tangent bundle
of $Y$. The representations of the covering $\com^*$
may be described by fractional weights of the original
$\com^*$.

The nonsingularity of $Y$ and $Y^f$ implies that there is a 
pullback from $A^{\com^*}_*(Y)$ to $A^{\com^*}_*(Y^f)$ and that,
 on each component
of $Y^f= \cup Y_i$,    
the usual self-intersection formula
$\iota^* \iota_* \alpha = e(N_i) \cdot \alpha$ holds.

Suppose  $\alpha=\sum \alpha_i$ 
in $A^{\C^*}_*(Y^f)_t$ 
pushes forward to zero.  Then, 
$$0 = \iota^* \iota_* \alpha = \sum e(N_i) \cdot \alpha_i$$
since pushing forward from one component of $Y^f$
and restricting to another necessarily gives zero.
Hence, for each $i$, $e(N_i) \cdot \alpha_i=0$.
However, since $e(N_i)$ is invertible in the localized
ring, each $\alpha_i$ is zero.
We have proven:

\begin{pr}
\label{nsloc}
The map $\iota_*: A^{\com^*}_*(Y^f) \rarr A^{\com^*}_*(Y)$ is an
isomorphism after localization.
\end{pr}

The self intersection formula also 
quickly implies the explicit localization
formula:
$$ [Y] = \sum_i \frac{[Y_i]}{e(N_i)}.$$
By Proposition \ref{nsloc},
there exists a unique class $\alpha$ satisfying $\iota_*(\alpha)=[Y]$. 
The condition $\iota^* \iota_*(\alpha)= [Y^f]$ then
determines $\alpha$.

If $X\subset Y$ is a $\com^*$-equivariant embedding,
the fixed substack $X^f$ may be {\em defined} by
$X^f= X \cap Y^f$. It follows from this
definition that
$$\Omega_X |^f_{X^f}= \Omega_{X^f}.$$ 
It is not difficult to show the
substack structure $X^f$ is independent of the
choice of nonsingular equivariant embedding. The constructions
and arguments for the virtual localization formula
for equivariant perfect obstruction theories on $X$
now go through unchanged.

\noindent

$$ $$

\noindent Department of Mathematics \\
University of Chicago \\
5734 S. University Ave. \\
Chicago, IL 60637 \\

\noindent graber@math.uchicago.edu \\
rahul@math.uchicago.edu

\end{document}